\documentclass[letterpaper, 10pt, conference]{ieeeconf}  % Comment this line out

\usepackage{cite}
\usepackage{graphicx}
\usepackage{url}
\usepackage{subcaption}
\usepackage{amsmath}
\usepackage{amsfonts}

\renewcommand{\baselinestretch}{1.1} 
\IEEEoverridecommandlockouts                              % This command is only
\newtheorem{assumption}{Assumption}

\title{\LARGE \bf Optimal Vehicle Dynamics and Powertrain Control for Connected and Automated Vehicles}

\author{Liuhui Zhao, {\itshape{Member, IEEE}}, A M Ishtiaque Mahbub, {\itshape{Student Member, IEEE}}, 
\\
Andreas A. Malikopoulos, {\itshape{Senior Member, IEEE}}
\thanks{This research was supported by ARPAE's NEXTCAR program under the award number DE-AR0000796.}% <-this % stops a space
\thanks{The authors are with the Department of Mechanical Engineering, University of Delaware, Newark, DE 19716 USA (emails: \tt\small{lhzhao@udel.edu; mahbub@udel.edu; andreas@udel.edu).}}%
}

\begin{document}

\maketitle
\thispagestyle{empty}

\begin{abstract}
The implementation of connected and automated vehicle technologies enables opportunities for a novel computational framework for real-time control actions aimed at optimizing energy consumption and associated benefits. In this paper, we present a two-level control architecture for a connected and automated
plug-in hybrid electric vehicle to optimize simultaneously its speed profile and powertrain efficiency. We evaluate the proposed architecture through simulation in a network
of vehicles. 
\end{abstract}

\indent

\section{Introduction}
In this paper, we are interested in investigating the opportunities to improve the efficiency of hybrid electric vehicles (HEVs) and plug-in HEVs (PHEVs) when these vehicles are connected and automated.
In an earlier work, we discussed the potential benefits of optimally coordinated connected and automated vehicles (CAVs) in a corridor using a traffic microsimulation environment without considering powertrain optimization \cite{Zhao2018}. In this paper, we apply a two-level supervisory control architecture for connected and automated PHEVs (CA-PHEVs) that consists of a vehicle dynamics (VD) controller and a powertrain (PT) controller. The supervisory controller oversees the VD and PT controllers and communicates the endogenous and exogenous information appropriately. The VD controller optimizes online the vehicle acceleration/deceleration profile to avoid stop-and-go driving in situations where there is a potential conflict with other vehicles, e.g., ramps, intersections, stop signs, roundabouts, etc. The PT controller computes the optimal nominal operation (set-points) for the engine and motor corresponding to the optimal solution derived from the VD controller. The complexity of the problem dimensionality can be managed by establishing two parallel and appropriately interacting computational levels, namely a cloud-based level, and a vehicle-based level. 

The objectives of this paper are to (1) optimize vehicle speed profile in terms of energy consumption and compare two different control approaches, and (2) optimize the power management of the vehicle for the optimal speed profile obtained in (1). We evaluate the effectiveness of the efficiency of the proposed architecture through simulation in a network of CA-PHEVs in MCity, a 32-acre vehicle testing facility located at the north campus of University of Michigan. The contribution of this paper is the analysis of online optimization of the vehicle- and powertrain-level operation of a CA-PHEV and classification of the improvements.

The remainder of the paper is organized as follows. In Section II, we summarize the relevant research efforts on vehicle dynamics control for CAVs and powertrain optimization for HEVs and PHEVs. In Section III, we introduce the control architecture that consists of the VD controller and the PT controller. In Section IV, we evaluate the effectiveness of the efficiency of the proposed approach in a simulation environment. Finally, we draw conclusions and discuss next steps in Section V.

\section{Literature Review}
%\subsection{Vehicle Dynamics Control}
Several research efforts have been reported in the literature proposing either centralized or decentralized approaches on coordinating CAVs. 
One of the very early efforts in this direction was proposed in 1969 by Athans \cite{Athans1969} for safe and efficient coordination of merging maneuvers with the intention to avoid congestion. %Assuming a given merging sequence, Athans formulated the merging problem as a linear optimal regulator to control a single string of vehicles, with the aim of minimizing the speed errors that will affect the desired headway between each consecutive pair of vehicles. 
Varaiya \cite{Varaiya1993} discussed extensively the key features of an automated intelligent vehicle-highway system and proposed a related control system architecture. Dresner and Stone \cite{Dresner2004} proposed the use of the reservation scheme to control a signal-free intersection. Other research efforts have focused on coordinating vehicles at intersections to improve the traffic flow \cite{Yan2009,kim2014}. More recently, online optimal control approaches were presented for coordinating CAVs in different transportation segments \cite{Rios-Torres2015, Ntousakis:2016aa, Rios-Torres2}. Closed-form, analytical solutions have been presented and tested in different transportation scenarios for coordinating online CAVs \cite{Zhang2016a,Malikopoulos2017, Malikopoulos2018a,Malikopoulos2018b, Malikopoulos2018c, Zhao2018, Mahbub2019CDC}. A survey paper \cite{Malikopoulos2016a} includes detailed discussions of the research efforts reported in the literature to date in this area.

%\subsection{Powertrain Optimization}
The development and implementation of a power management control algorithm online constitutes a challenging control problem and has been the subject of intense study for the last two decades \cite{Malikopoulos2016AMO, Malikopoulos2014c,Malikopoulos2013a, Malikopoulos2015_ITS_HEV}. 
%\cite{Malikopoulos2013a,Malikopoulos2016AMO,Malikopoulos2014c,Malikopoulos2014dec,Malikopoulos2015_ITS_HEV}
Several research efforts have focused on different heuristic approaches to optimized the power management control in PHEVs \cite{Saeks2002, Schouten2002}. The latter, however, cannot completely encompass all the potential fuel economy and emission benefits of the various HEV architectures. 
%Dynamic programming (DP), on the other hand, can yield global optimal solution for both stochastic and deterministic formulation, and provide a rigorous approach for sequential decision-making problems under uncertainty. 
A significant amount of work on optimizing the power management control in HEVs has focused on applying dynamic programing (DP) \cite{Lin2003,Tate2010} to benchmark fuel economy. To address the computational constraints associated with DP, research efforts have been concentrated on developing online algorithms consisting of an instantaneous optimization problem. Over the years, different online control approaches based on Pontryagin's Maximum Principle  \cite{Serrao2008}, Equivalent Consumption Minimization Strategy (ECMS) \cite{Paganelli2001, Sciarretta2004}, and Adaptive-ECMS \cite{Musardo2005} have been proposed. Recently, there have been efforts to address the vehicle-level and powertrain-level operation simultaneously \cite{Luo2015b, Li2017a}. However, these efforts have exhibited some limitations for online implementation. 
A detailed discussion of the research reported in the literature today on the power management control for HEVs/PHEVs can be found in \cite{Malikopoulos2014b}.

\section{Control Architecture}
\begin{figure}[ht]
\centering
\includegraphics[width=0.38\textwidth]{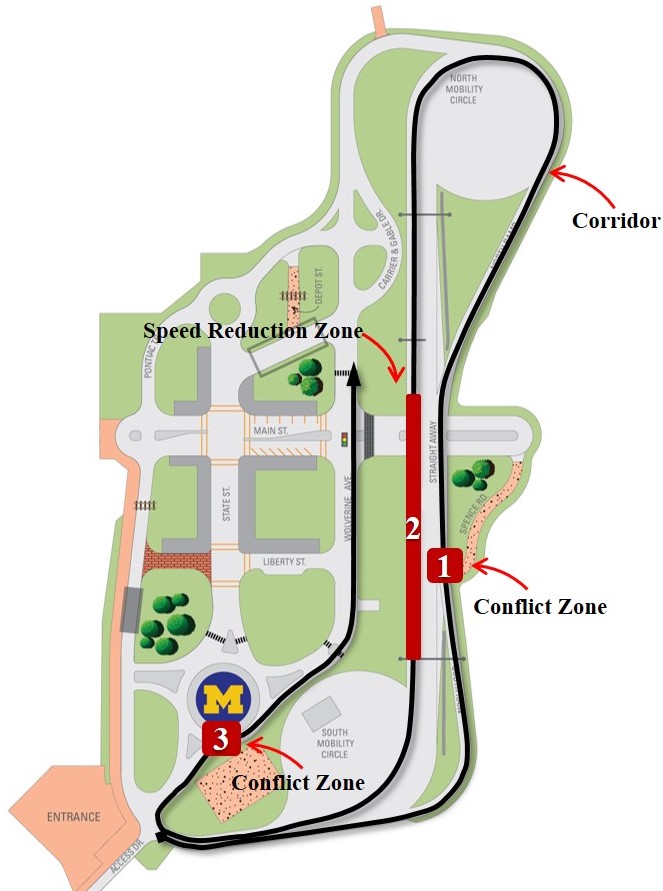} 
\caption{The corridor in Mcity with the conflict zones.}%
\label{fig:corridor}%
\end{figure}

We consider a network of CA-PHEVs cruising through a corridor in Mcity that consists of several conflict zones, e.g., a highway on-ramp, a speed reduction zone (SRZ), and a roundabout as shown in Fig. \ref{fig:corridor}. The CA-PHEVs are retrofitted with necessary communication devices to interact with other CA-PHEVs and infrastructure within their communication range. The proposed control architecture applies to the operation of each CA-PHEV over a range of real-world driving scenarios deemed characteristic of typical commutes. 
%A typical commute of a vehicle includes merging at roadways, crossing signalized intersections, cruising in congested traffic, and passing through speed reduction zones. Therefore, t
The vehicle dynamics (VD) controller optimizes the speed profile of the vehicle over such traffic scenarios. The powertrain (PT) controller computes the optimal nominal operation (set-points) for the engine and motor corresponding to the optimal speed profile provided by the VD controller. 
%The supervisory controller coordinates the VD and PT controllers to ensure the optimal solution yielded by the VD controller is feasible for the PT controller. However, the details of this coordination along with the implications of the computational efforts of the VD and PT controllers are outside the scope of this paper.

\subsection{Vehicle Dynamics (VD) Controller}
\subsubsection{Vehicle Dynamics Model}
Let $\mathcal{N}(t)={1,\dots,N(t)}$, where $N(t)$ is the total number of vehicles inside the control zone at time $ t\in \mathbb{R}^+$, be a queue of vehicles to be analyzed. The dynamics of each vehicle $ i\in \mathcal{N}(t)$ are represented with a state equation
\begin{equation} \label{eq:state}
\dot{x}_i(t) = f(t, x_i, u_i), ~ x_i(t_i^0) = x_i^0,
\end{equation}
where $x_i(t), u_i(t)$ are the state of the vehicle and control input, $t_i^0$ is the initial time of vehicle $i\in\mathcal{N}(t)$, and $x_i^0$ is the value of the initial state. For simplicity, we model each vehicle as a double integrator, i.e., $\dot{p}_i = v_i(t)$ and $\dot{v}_i = u_i(t)$, where $p_i(t) \in \mathcal{P}_i, v_i(t) \in \mathcal{V}_i$, and $u_i(t) \in \mathcal{U}_i$ denote the position, speed, and acceleration/deceleration (control input) of each vehicle $i$. Let $x_i(t)=[p_i(t)~ v_i(t)]^T$ denotes the state of each vehicle $i$, with initial value $x_i^0=[0 ~ v_i^0]^T$, taking values in the state space $\mathcal{X}_i=\mathcal{P}_i\times  \mathcal{V}_i$. The sets $\mathcal{P}_i, \mathcal{V}_i$, and $\mathcal{U}_i, i\in \mathcal{N}(t)$, are complete and totally bounded subsets of $\mathbb{R}$. The state space $\mathcal{X}_i$ for each vehicle $i$ is closed with respect to the induced topology on $\mathcal{P}_i\times  \mathcal{V}_i$ and thus, it is compact.

To ensure that the control input and vehicle speed are within a given admissible range, the following constraints are imposed.
\begin{equation}\label{eq:constraints}
\begin{aligned} 
u_{min} &\leq u_i(t) \leq u_{max},~ \text{and} \\
0 &\leq v_{min} \leq v_i(t) \leq v_{max}, ~ \forall t \in [t_i^0, ~ t_i^f]
\end{aligned} 
\end{equation}
where $u_{min}, u_{max}$ are the minimum deceleration and maximum acceleration respectively, $v_{min}, v_{max}$ are the minimum and maximum speed limits respectively, and $t_i^0$, $t_i^f$ are the times that each vehicle $i$ enters and exits the corridor.

To ensure the absence of rear-end collision of two consecutive vehicles traveling on the same lane, the position of the preceding vehicle should be greater than, or equal to the position of the following vehicle plus a predefined safe distance $\delta_i(t)$, where $\delta_i(t)$ is proportional to the speed of vehicle $i$, $v_i(t)$. Thus, we impose the rear-end safety constraint
\begin{equation} \label{eq:safety}
s_i(t)=p_k(t)-p_i(t)\geq \delta_i(t), \forall t \in  [t_i^0, ~ t_i^f]
\end{equation}
where vehicle $k$ is immediately ahead of $i$ on the same lane. 

In the modeling framework described above, we impose the following assumptions:

\begin{assumption}
For each CA-PHEV $i$, none of the constraints is active at $t_i^0.$
\end{assumption}

\begin{assumption}
Each vehicle is equipped with sensors to measure and share their local information while communication among CA-PHEVs occurs without any delays or errors. 
\end{assumption}

The first assumption ensures that the initial state and control input are feasible. The second assumption might be strong, but it is relatively straightforward to relax as long as the noise in the measurements and/or delays is bounded. 
%For example, we can determine upper bounds on the state uncertainties as a result of sensing or communication errors and delays, and incorporate these into more conservative safety constraints.

Next, we present two approaches for the VD controller, using isolated conflict zones and a coordinated corridor control approach. In the former approach, we assume that the conflict zones in the corridor are independent, such that vehicles only coordinate with each other to travel through an immediate downstream conflict zone. In the latter approach, we assume that all the vehicles traveling in the corridor coordinated with each other to pass through all the conflict zones smoothly.

\subsubsection{Isolated Conflict Zones}
We consider a corridor that contains three conflict zones (Fig. \ref{fig:corridor}), e.g., a merging roadway (conflict zone 1), a SRZ (conflict zone 2), and a roundabout (conflict zone 3). Upstream of each conflict zone, we define a \textit{control zone} where CA-PHEVs coordinate with each other to avoid any rear-end or lateral collision in the conflict zone. For each conflict zone, there is a coordinator that communicates with the CA-PHEVs traveling within the control range. Note that the coordinator is not involved in any decision of the vehicles. The coordinator just assigns a unique identity to each CA-PHEV when they enter a control zone. The objective of each CA-PHEV is to derive its optimal control input (acceleration/deceleration) to cross the conflict zones without any rear-end or lateral collision with the other vehicles. 

Let $z \in \mathbb{Z}$ be the index of a conflict zone in the corridor. When a vehicle enters the control zone, the coordinator receives its information and assigns a unique identity $i$ to the vehicle. Let $t_i^{z, 0}$ be the time when vehicle $i\in\mathcal{N}(t)$ enters the control zone towards conflict zone $z$, and $t_i^{z,f}$ be the time when vehicle $i$ exits the corresponding control zone. In each control zone, we denote the sequence of the vehicles to be entering a conflict zone as $\mathcal{N}_z(t)={1,..., N(t)}$. Thus, we formulate the following optimization problem for each vehicle in the queue $\mathcal{N}_z(t)$
\begin{gather} \label{eq:min}
\min_{u_i}\frac{1}{2}\int_{t_i^{z, 0}}^{t_i^{z,f}} u_i^2(t)dt, ~\forall i \in \mathcal{N}_z(t), ~\forall z \in \mathcal{Z} \\
\text{subject to}: (\ref{eq:state}), (\ref{eq:constraints}), (\ref{eq:safety}), \nonumber\\
p_{i}(t_i^{z,0})=p_{i}^{z,0}\text{, }v_{i}(t_i^{z,0})=v_{i}^{z,0}\text{, }p_{i}(t_i^{z,f})=p_{z},\nonumber\\
\text{and given }t_i^{z,0}\text{, }t_i^{z,f},\nonumber
\end{gather}
where $p_{z}$ is the location (i.e., entry position) of conflict zone $z$, $p_{i}^{z,0}$, $v_{i}^{z,0}$ are the initial position and speed of vehicle $i$ when it enters the control zone of conflict zone $z$. 
To address this problem, we apply the decentralized optimal control framework and analytical solution presented in \cite{Rios-Torres2,Malikopoulos2017,Malikopoulos2018c}. %The solution of the constrained problem has been presented in \cite{Malikopoulos2017}. %The problem of minimizing the control input (acceleration/deceleration) for each vehicle $i$ from the time $t_i^{z, 0}$ that the vehicle enters the control zone until the time $t_i^{z, f}$ that it exits the control zone. 

\subsubsection{Coordinated Corridor Control Approach}
In this approach, we consider a single coordinator that monitors all  vehicles traveling along the corridor. Note that the coordinator serves as an information center which is able to collect vehicular data through vehicle-to-infrastructure communication and is not involved in any decision on the vehicle operation. Thus, the coverage of the coordinator is flexible and the length of corridor could be extended in the presence of connected infrastructure. %The details of upper level scheduling policy and  the lower level optimal control problem for every vehicle, i.e., the optimal acceleration profile that will achieve the pre-determined entry time at each conflict zone, was presented in \cite{Zhao2018}.

Let $\mathcal{N}(t)\in\mathbb{N}$ be the number of CA-PHEVs in the corridor at time $t\in\mathbb{R}^{+}$. When a vehicle enters the boundary of the corridor, it broadcasts its route information to the coordinator. Then, the coordinator assigns a unique identification number $i\in\mathbb{N}$. Let $t_i^0$ be the initial time that vehicle $i$ enters the corridor, $t_{i}^{z}$ be the time for vehicle $i$ to enter the conflict zone $z$, $z \in \mathcal{Z}$, and $t_{i}^{f}$ be the time for vehicle $i$ to enter the final conflict zone.

To avoid any possible lateral collision, there is a number of ways to compute $t_{i}^{z}$ for each CAV $i$. In this paper, we adopt a scheme reported in \cite{Zhao2018} for determining the time $t_{i}^{z}$ that each CAV enters the conflict zone $z \in \mathcal{Z}$. Thus, for each vehicle, we formulate the following optimal control problem the solution of which  yields for each vehicle the optimal control input (acceleration/deceleration) to achieve the assigned time $t_{i}^{z}$ (upon arrival of CAV $i$) without collision
\begin{gather} \label{eq:min_cor}
\min_{u_i}\frac{1}{2}\int_{t_i^0}^{t_i^f} u_i^2(t)dt, ~\forall i \in \mathcal{N}(t) \\
\text{Subject to}: (\ref{eq:state}), (\ref{eq:constraints}),  (\ref{eq:safety}), \nonumber\\
p_{i}(t_i^0)=p_{i}^{0}\text{, }v_{i}(t_i^{0})=v_{i}^{0}\text{, }p_{i}(t_i^{z})=p_{z},\nonumber\\
\text{and given }t_i^{0}\text{, }t_i^{z}\text{, }t_i^{f}, ~\forall z \in \mathcal{Z},\nonumber
\end{gather}
where $p_{i}^{0}$, $v_{i}^{0}$ are the initial position and speed of vehicle $i$ when it enters the corridor, and $t_i^{z}$ is the time when vehicle $i$ enters the conflict zone $z$. 

Each vehicle $i\in\mathbb{N}$ determines the time $t_i^{z}$ that will be entering the conflict zone $z \in \mathcal{Z}$ upon arrival at the entry of the corridor. Thus, the next vehicle $i+1$, upon its arrival at the entry of the corridor, will search for feasible times to cross the conflict zones based on available time slots. The recursion is initialized when the first vehicle enters the control zone, i.e., it is assigned $i = 1$. The details of determining the sequence that each vehicle enters the control zone along with the closed-form analytical solution can be found in \cite{Zhao2018}.

\subsection{Powertrain (PT)  Controller}
The objective of the PT controller is to derive an optimal control policy to split the torque demanded by the driver, $T_{driver}$ between the engine and the motor torque, $T_{eng}$ and $T_{mot}$ respectively, for the optimal speed profile derived by the VD controller.
We model the evolution of the CA-PHEV state at each stage $t=0,1,...$ as a controlled Markov chain following the framework reported in \cite{Malikopoulos2016AMO}, with a finite state and control space,  $\mathcal{S} \subset \mathbb{R}^n$, and  $\mathcal{U} \subset \mathbb{R}^m, n,m\in\mathbb{N}$ respectively. We introduce the sequence of the random variables  $X_{t(1:2)}=(X_{t(1)},  X_{t(2)}))^T=(N_{eng},  N_{mot})^T\in\mathcal{S}$, $U_{t(1:2)}=(U_{t(1)},  U_{t(2)}))^T=(T_{eng},  T_{mot})^T\in\mathcal{U}$ and $W_{t(1:2)}$, corresponding to the HEV state (engine and motor speed, $N_{eng}$ and $N_{mot}$), control action (engine and motor torque, $T_{eng}^*$ and $T_{mot}^*$) and system uncertainty in terms of the torque demanded by the driver as designated by the pedal position respectively. At each stage $t$, the controller observes the system state  $X_{t(1:2)}=i\in\mathcal{S}$, and executes an action, $U_{t(1:2)}\in\mathcal{U}$ at that state. At the next stage, $t+1$, the system transits to the state $X_{t+1(1:2)} =j\in\mathcal{S}$ and a one-stage expected cost, $k(X_{t(1:2)},U_{t(1:2)})$, is incurred corresponding to the engine's and motor's efficiency. After the transition to the next state, a new action is selected and the process is repeated. 

We select the average cost criterion as we wish to optimize the efficiency of each CA-PHEV on average. Hence, the objective of the PT controller to derive a stationary control policy that minimizes the long-run expected average cost,
\begin{equation}
J^\pi= \lim_{T\to\infty}\frac{1}{T+1}\mathbb{E}^\pi \left [\sum_{t=0}^{T} k(X_{t(1:2)},U_{t(1:2)})\right],\label{eq:2}
\end{equation}
where $k(X_{t(1:2)},U_{t(1:2)})$ is the one-stage cost of CA-PHEV. However, onboard derivation of the control policy by minimizing \eqref{eq:2} is not feasible due to the associated computational burden. Therefore, we formulate a multiobjective optimization problem that yields the Pareto control policy. It has be proven that the Pareto policy is equivalent to the optimal control policy derived by solving \eqref{eq:2} with DP \cite{Malikopoulos2016b}.
To this end, we formulate the multiobjective optimization problem by considering the engine's efficiency, $\eta_{eng},$ and the motor's efficiency, $\eta_{mot}$. As the engine's efficiency is a function of engine torque $T_{eng}$ and engine speed $N_{eng}$, we construct the engine's efficiency as $f_1(N_{eng}, T_{eng}) = \eta_{eng}$. Due to similar dependency, we also write the motor's efficiency as, $f_2(N_{mot},T_{mot}) = \eta_{mot}$.

The multiobjective optimization problem is formulated as
\begin{gather} 
\min_{U_t} k(X_{t(1:2)}, U_{t(1:2)})= \nonumber\\
\max_{U_t} \big(\alpha\cdot f_1(X_{t(1)}, U_{t(1)}) + (1-\alpha)\cdot f_2(X_{t(2)}, U_{t(2)})\big)\nonumber \\
\text{s.t.} \ \sum_{i=1}^{2}U_{t(i)}= T_{driver},\label{eq:objective_function}
\end{gather}
where $\alpha$ is a scalar that takes values in [0,1], $f_1(N_{eng}, T_{eng}) = \eta_{eng}$, $f_2(N_{mot},T_{mot}) = \eta_{mot}$ are the engine's efficiency and  motor's efficiency constituting the multiobjective function, and $T_{driver}$ is the torque demanded by the driver. The solution of the multiobjective optimization problem in \eqref{eq:objective_function} yields the Pareto efficiency set between the engine and the motor by varying $\alpha$ from $0$ to $1$ at any given state of the CA-PHEV and torque demand. We calculate the Pareto efficiency set offline and store it in a table onboard the vehicle. Then, we implement the optimal policy online using this table. 

\section{Simulation Evaluation}
%
%\begin{figure}[hb]
%	\centering
%	\includegraphics[width=0.45\textwidth]{figs/pareto_engine}
%	\caption{Pareto efficiency set of the powertrain controller.}%
%	\label{fig:pareto_result}%
%\end{figure}

\subsection{Simulation Environment}

To evaluate the optimal control policy of the PT controller we use a modified version of the Vehicle-Engine SIMulation (VESIM) model reported in \cite{Malikopoulos:2006aa} and references therein.
We primarily focus on two objectives: 1) to evaluate the network performance with the proposed supervisory control framework and 2) to evaluate the efficiency of different VD controller approaches under three different traffic levels (i.e., light, medium, and heavy traffic conditions). To this end, we develop three scenarios. For each scenario, we apply both baseline (case with human-driven vehicles) and optimal PT controller to estimate fuel efficiency for all CA-PHEVs traveling in the corridor. The corridor has a length of 1.3 $km$ in MCity (Fig. \ref{fig:corridor}). The desired speeds for the highway, urban and SRZ are set as 17 $m/s$, 11 $m/s,$ and 8 $m/s$ respectively. The length of the speed reduction zone is 125 $m$. 

\begin{itemize}
\item{Scenario 1 (Baseline): All vehicles in the network are human-driven vehicles. In this scenario, the Wiedemann car following model \cite{wiedemann1974} built in VISSIM is applied. The intersection is controlled by fixed-time signal controller, whose signal timing is optimized for the traffic condition set in the study.}
\item{Scenario 2 (Isolated conflict zones): We consider \\100\% penetration rate of CA-PHEVs in this scenario. The vehicles optimize their trajectory based on the isolated conflict zone control approach presented in Section III.}
\item{Scenario 3 (Coordinated corridor control approach): We consider 100\% penetration rate of CA-PHEVs in this scenario.  The vehicles optimize their trajectory based on the coordinated corridor control approach presented in Section III.}
\end{itemize}

Under Scenario 2, the length of the control zone for the on-ramp, the SRZ and the roundabout were selected to be 150 $m$. 

\subsection{Evaluation Results}
\subsubsection{The effectiveness of the VD Controller} 
\begin{figure*}[t]
	\centering
	\begin{subfigure}[b]{0.29\textwidth}
		\includegraphics[width=\textwidth]{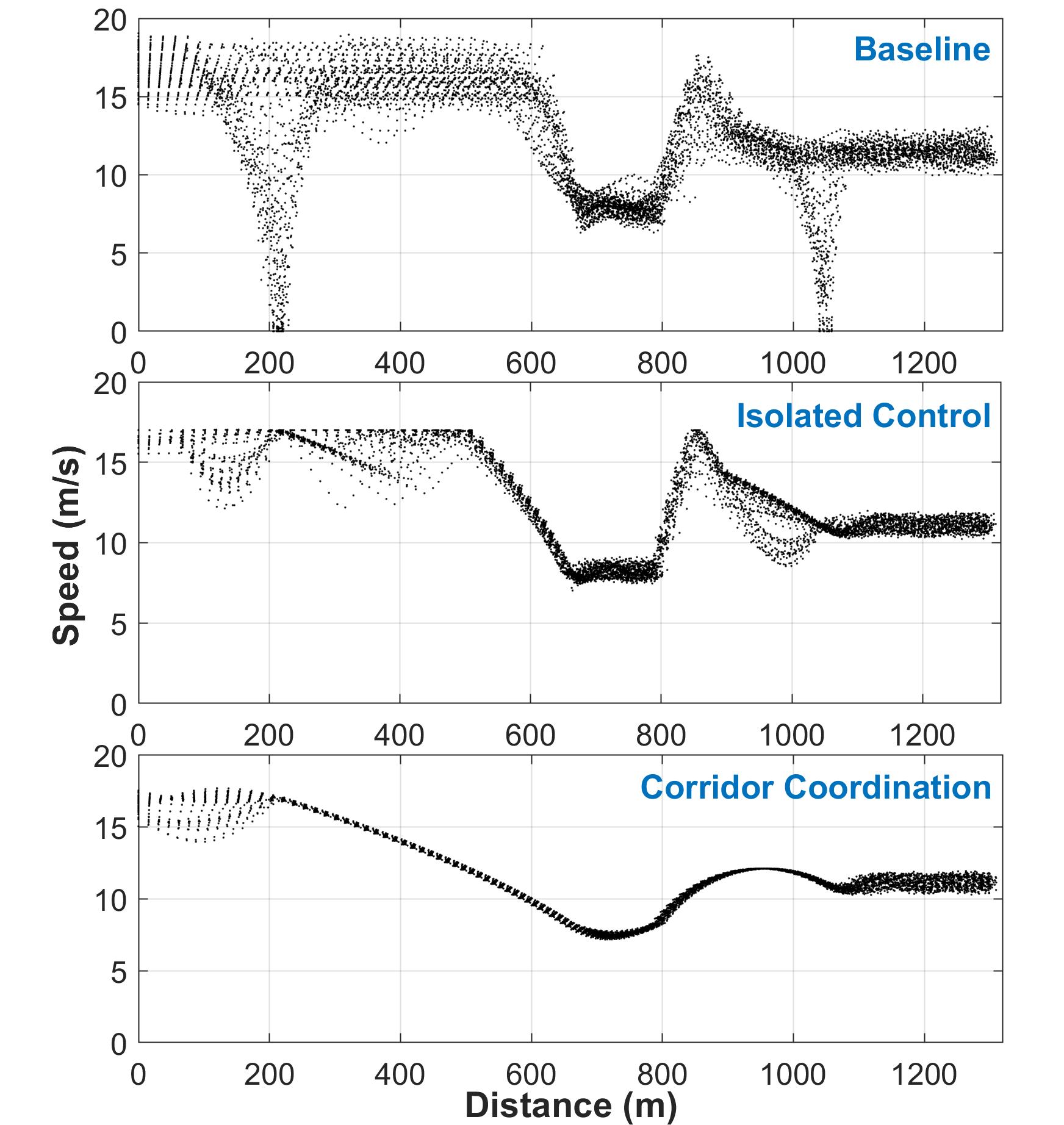} 
		\caption{light traffic demand level} \label{fig:light}
	\end{subfigure}	
	\begin{subfigure}[b]{0.29\textwidth}
		\includegraphics[width=\textwidth]{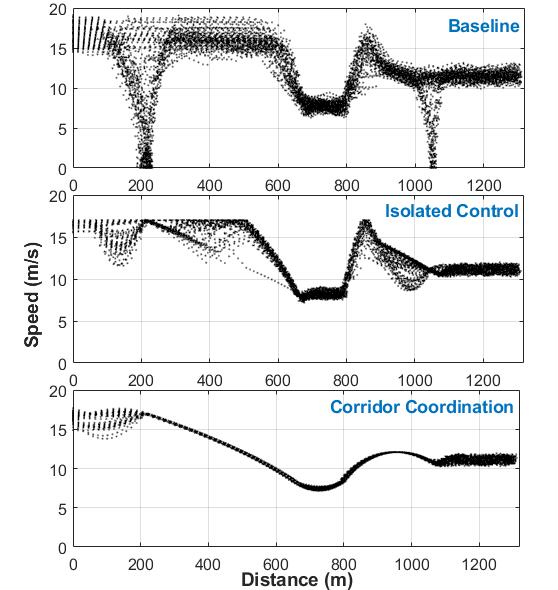} 
		\caption{medium traffic demand level} \label{fig:med}
	\end{subfigure}	
	\begin{subfigure}[b]{0.29\textwidth}
		\includegraphics[width=\textwidth]{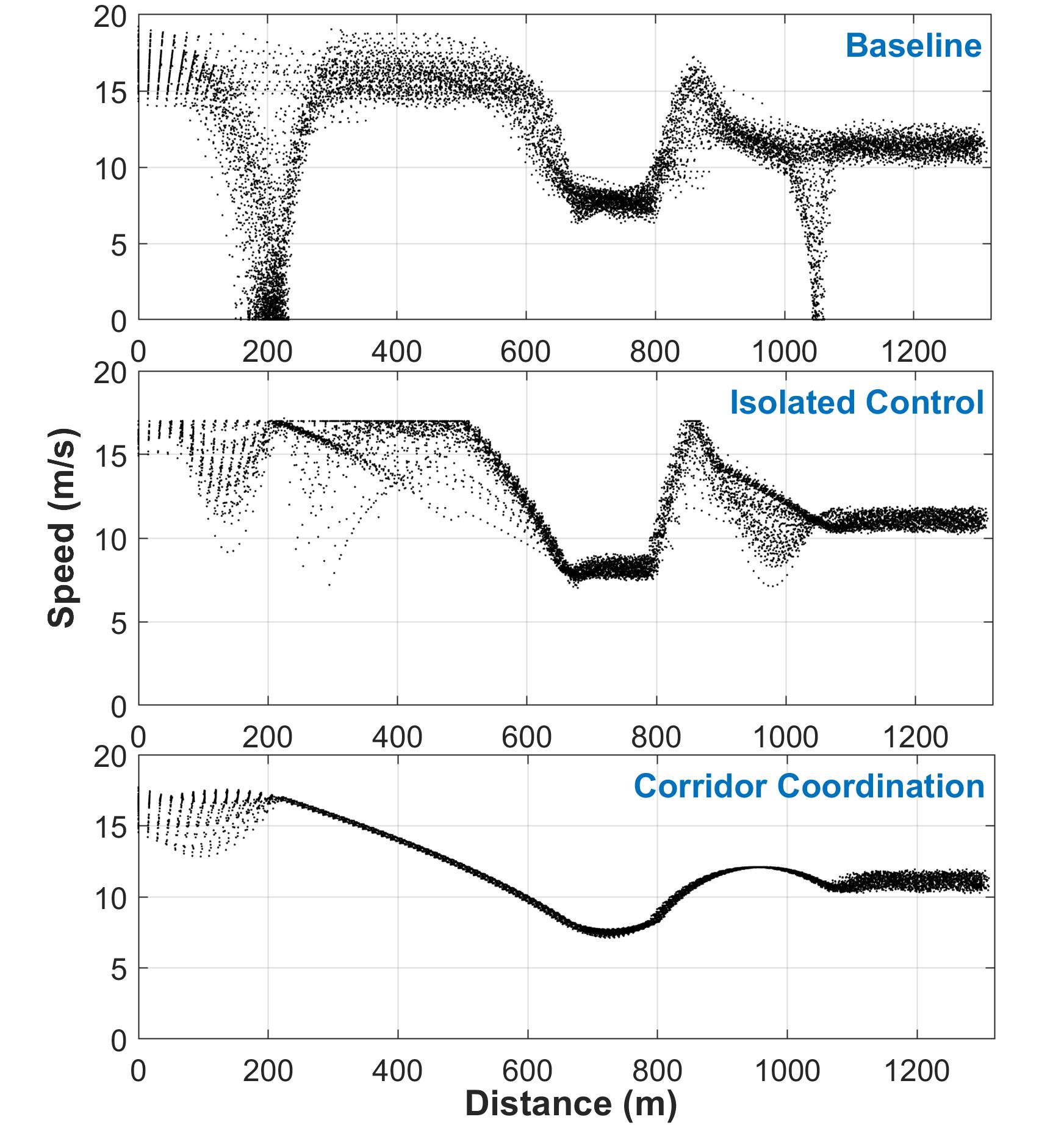} 
		\caption{heavy traffic demand level} \label{fig:high}
	\end{subfigure}
	
	\caption{Speed profiles for all CA-PHEVs under different VD control approaches and demand levels.}
	\label{fig:speed}%
\end{figure*}

The vehicle trajectories are shown in Fig. \ref{fig:speed}. Under Scenario 1, the vehicles need to slow down or stop for merging into highway or roundabout. Therefore, we note (top panel of Fig. \ref{fig:speed}) that there are many stop-and-go events in vehicle speed profiles under the baseline scenario. 
Under scenario 2 (middle panel of Fig. \ref{fig:speed}), we observe smooth speed profiles inside each control zone. At the downstream of each conflict zone, vehicles exit from the VD control zone, thus we see similar speed patterns as under scenario 1 outside the control zones. We also note that while the isolated control approach is able to eliminate stop-and-go driving, the resulting increased traffic flow into the downstream speed reduction zone leads to speed reduction downstream the highway on-ramp segment.

On the contrary, with the coordinated control scheme under scenario 3, traffic information of the entire corridor is shared among all vehicles. Therefore, all the CA-PHEVs travel through the corridor are able to optimize their trajectories at the entry of the corridor and drive smoothly throughout the corridor, even with high traffic demand level (as shown in the bottom panel of Fig. \ref{fig:speed}). 

%Under three different traffic conditions, we calculate vehicle travel times through the corridor under different VD control approaches and plot travel time distribution in Fig.~\ref{fig:time}. We observe that with coordinated corridor control approach, due to longer preparation for downstream roadway segment, the average travel time under light traffic condition is longer than the baseline scenario. However, the results reveal that through the corridor coordination, the traffic flow is further smoothed (much lower variation in vehicle travel times as shown in Fig.~\ref{fig:time}) compared to the isolated conflict zone control approach. 

%\begin{figure*}[h]
%\centering
%\includegraphics[width=.95\textwidth]{figs/travel_time_overall}
%\caption{Travel time distribution under different VD control approaches and traffic demand levels.}
%\label{fig:time}
%\end{figure*}

%\begin{figure}[!h]
%\centering
%\includegraphics[width=0.45\textwidth]{figs/veh149_spd}
%\caption{Vehicle speed profiles under different VD control approaches.}
%\label{fig:veh_spd}
%\end{figure}

\subsubsection{The effectiveness of PT Controller} 

\begin{table}[b]
\centering
\caption{Improvements of fuel efficiency over baseline scenarios under different traffic conditions.}
\includegraphics[width=0.48\textwidth]{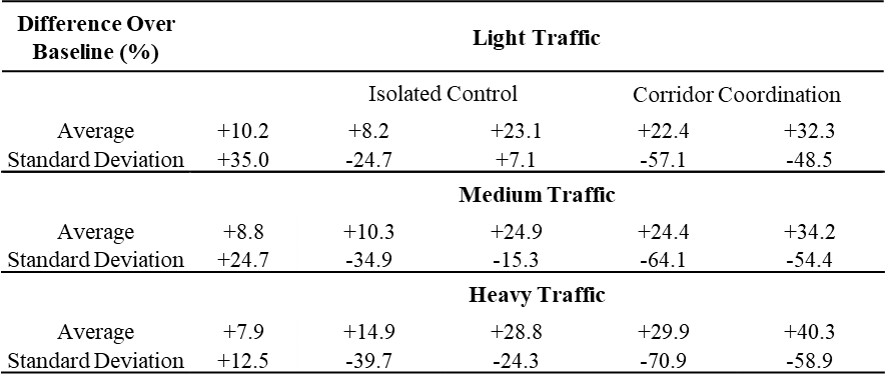}
\label{tab:fuel_comp}
\end{table}

Table \ref{tab:fuel_comp} summarizes the results for the baseline and optimal PT controllers for all VD control approaches under different traffic conditions. Without the PT controller, the coordinated corridor control approach leads to significantly efficiency improvement in comparison with the isolated control approach. This improvement can be attributed to the fact that the coordinated corridor control approach eliminates all the uncontrolled patches between the control zones of individual coordinators in isolated control approach and optimizes the entire trajectory inside the corridor, leading to improved fuel efficiency. At the same time, as all the vehicles fall under a single coordinated optimal control policy inside the corridor, significant decrease in terms of the standard deviation of average efficiency is observed as compared to the isolated control case (Table \ref{tab:fuel_comp}). 
We observe that for the optimal speed profiles yielded from the corridor coordination control approach, the contribution of optimal PT controller in terms of overall fuel efficiency improvement is lower than those yielded from the isolated control approach. The reason is that under the corridor coordination approach, each vehicle spend a significant amount of time in the regenerative braking event, which reduces opportunity for the Pareto optimal policy of the  PT controller to be active. Note that, the optimal PT controller coupled with different VD control policy shows consistent improvement in terms of efficiency in all traffic conditions for isolated control and for corridor coordination policy. As both the VD control policy eliminate stop-and-go driving and smooths out the velocity profile, the improvement in average efficiency remains consistent irrespective of traffic volumes.

\section{Concluding Remarks}
In this paper, we presented a two-level control architecture with the aim to optimize simultaneously vehicle-level and powertrain-level operation of a CA-PHEV.
We used two different control approaches to optimize simultaneously the vehicle speed profile and the powertrain efficiency of a CA-PHEV. 
%We also compared the efficiency of different VD control approaches, and showed that although the coordinated corridor approach improves significantly fuel efficiency, vehicle travel time savings are marginal compared to the isolated conflict zone control approach.
We applied the proposed architecture to the operation of CA-PHEVs over a range of real-world driving scenarios deemed characteristic of typical commute that included a merging roadway, a speed reduction zone, and a roundabout. 
Ongoing work considers additional scenarios, including intersections, while incorporating the state and control constraints in the analytical solution of the VD controller. 
Although potential benefits of full penetration of CA-PHEVs to alleviate traffic congestion and reduce fuel consumption have become apparent, different penetrations of CA-PHEVs can alter significantly the efficiency of the entire system. Future work should consider the interactions between CA-PHEVs and human-driven vehicles using feedback for the drivers \cite{Malikopoulos2012a}, and the implications of different  penetration of CA-PHEVs.  Although it is relatively straightforward to extend our results to the case where the perfect information assumption is relaxed, future research needs to be directed at the implications of errors and/or delays.

\bibliographystyle{IEEEtran}
\bibliography{ccta_corridor_ref}

\end{document}